\documentclass[a4paper]{article}
\usepackage[comma,numbers,square,sort&compress]{natbib}

\newtheorem{theorem}{Theorem}

\newtheorem{remark}{Remark}

\newtheorem{corollary}{Corollary}

\renewcommand{\top}{{\mathrm{T}}}
\usepackage{graphicx,amsmath,amssymb}
\usepackage{subfigure,pgf,tikz,pst-node}
\usetikzlibrary{arrows,shapes,automata,positioning,fit}

\def\rm{\mathrm}
\renewcommand{\top}{{\mathrm{T}}}

\begin{document}

\title{Output consensus of nonlinear multi-agent systems with unknown control directions}

\author{Yutao Tang \footnote{This work was partially supported by the Fundamental Research Funds for the Central Universities under Grant No.~24820152015RC36.  Yutao Tang is with the School of Automation, Beijing University of Posts and Telecommunications, Beijing 100876. P.\,R. China. Email: yttang@bupt.edu.cn}}

\date{}

\maketitle

{\noindent\bf Abstract}: In this paper, we consider an output consensus problem for a general class of
nonlinear multi-agent systems without a prior knowledge of the agents'
control directions. Two distributed Nussbaum-type control laws are proposed
to solve the leaderless and leader-following adaptive consensus for
heterogeneous multiple agents. Examples and simulations are given to verify
their effectiveness.


\section{Introduction}

Multi-agent system has been a hot topic in the last decades due to its
numerous applications, such as cooperative control of unmanned aerial
vehicles, communication among sensor networks, and formation of mobile
robots(\cite{olfati2007consensus,ren2008distributed}). As one of the most
important problems, consensus with or without leaders has been extensively
studied. For example, various distributed protocols were proposed and
analyzed in \cite{ren2008distributed} for both leader-following and
leaderless cases. In \cite{hong2006tracking}, the authors proposed a
distributed observer-based control law using local information to track an
integrator-type leader. Later, this work has been extended to multi-agent
systems with a general linear dynamics (\cite{hong2009multi, ni2010leader})
under switching topologies.  Other extensions including consensus under
time-delay communications or with an unknown-input driven leader have been
studied in \cite{gao2011consensus, hu2007leader, tang2014leader}.

Output consensus problem of nonlinear agents have also been studied by many authors. As an extension for single integrators, the output synchronization of a group of input-output passive nonlinear systems were investigated in \cite{chopra2006passivity} by adding proper couplings between them. Similar results were obtained in \cite{hamadeh2012global}, where networks of cyclic feedback biochemical oscillators were analyzed by dissipativity theory to achieve synchronization. In \cite{wang2014consensus}, nonlinear multi-agent systems in output feedback form with unity relative degree was studied to achieve leader-follower consensus. However, most of existing results were obtained by assuming the control direction matrix or at least the sign for single-input single-output case was known.

In practice, a control direction may not always be known a priori in many applications. For example, under some steering conditions like a course-changing operation, the control direction of a ship may be unknown \cite{du2007adaptive}. Even if it is known at first, the control direction of a plant may be changed by some structural damage \cite{liu2008multivariable}. To tackle this problem, the Nussbaum-type function, originally proposed in \cite{nussbaum1983some}, has been extensively used to solve such a problem.

For a single system with only one Nussbaum-type function used (\cite{nussbaum1983some}), the stability of the closed-loop system can be analyzed relatively easily. For a group of systems that are (physically) interconnected and each system involves an unknown control direction, the control design and stability analysis become more difficult since the controller of each subsystem might need one Nussbaum-type function.  A decentralized adaptive control problem was first tackled in \cite{ye1999decentralized} by respectively constructing a single Lyapunov function candidate for each subsystem and employing one Nussbaum-type function for each corresponding subsystem.  Recently, to achieve the consensus of integrator-type agents, a special type of Nussbaum functions has been constructed in \cite{chen2014adaptive}. However, these results relied on the assumption that all unknown control directions are with the same sign.  In fact, the results on consensus for agents with unknown heterogonous high-frequency gain signs are very few with only one exception \cite{peng2014cooperative}, where a Nussbaum-type adaptive controller was designed for each single-integrator agent such that consensus of the multi-agent network can be achieved. Nevertheless, the solvability of consensus among more general linear and nonlinear agents with unknown control directions is still unclear.

To answer this question, we mainly emphasize on a class of nonlinear heterogeneous agents having a passivity-like property with unknown control directions. By constructing a Nussbaum-type protocol, we can achieve leaderless output consensus under a mild connectivity condition of the communication graph among those agents.  With some minor modifications, this protocol is also able to drive all agents to some desired value, such as the equilibrium. In this way, our main contributions are at least two-fold:
\begin{itemize}
  \item We consider a group of heterogeneous agents with unknown control directions. Comparing with the results in \cite{chen2014adaptive}, this control protocol can allow the cases when those control directions have nonidentical signs. These conclusions are consistent with existing coordination results in \cite{bai2011cooperative,ren2008distributed}.
  \item We consider a more general class of agent's systems with unknown control directions, which strictly cover those heterogeneous single integrators considered in \cite{chen2014adaptive, peng2014cooperative}. Even for those single integrators, these results still hold for a class of directed graphs, while the conclusions in \cite{chen2014adaptive} are obtained only for undirected graphs.
\end{itemize}

The rest of this paper is organized as follows. Some preliminaries and problem formulation are given in Section 2 and 3. Main results are presented in Section 4, where two types of adaptive control laws are proposed for both leaderless and leader-following cases. Finally, simulations and our concluding remarks are presented at the end.

Notations: Let $\mathbb{R}^n$ be the $n$-dimensional Euclidean space, and $\mathbb{R}^{n\times m}$ be the set of $n\times m$ real matrices. For a vector $x$, $||x||$ denotes its Euclidian norm. $\text{diag}\{b_1,{\dots},b_n\}$ denotes an $n\times n$
diagonal matrix with diagonal elements $b_i\; (i=1,{\dots},n)$;
$\text{col}(a_1,{\dots},a_n) = [a_1^\top,{\dots},a_n^\top]^\top$ for any column
vectors $a_i\; (i=1,{\dots},n)$.

\section{Preliminaries}

Before the main results, we introduce some preliminaries on graph theory and nonlinear systems.

A directed graph (or digraph) $\mathcal {G}=(\mathcal {N}, \mathcal {E})$, where $\mathcal{N}=\{1,{\dots},n\}$ is the set of nodes and $\mathcal {E}$ is the set of edges (\cite{mesbahi2010graph}). $(i,j)$ denotes an edge leaving from node $i$ and entering node $j$. A directed path in graph $\mathcal {G}$ is an alternating sequence $i_{1}e_{1}i_{2}e_{2}{\cdots}e_{k-1}i_{k}$ of nodes $i_{l}$ and edges $e_{m}=(i_{m},i_{m+1}) \in\mathcal {E}$ for $l=1,2,{\dots},k$. If there exists a directed path from node $i$ to node $j$ then node $i$ is said to be reachable from node $j$. The neighbor set of agent $i$ is defined as $\mathcal{N}_i=\{j: (j,i)\in \mathcal {E} \}$ for $i=1,...,n$. A weighted adjacency matrix of a digraph $\mathcal {G}$ is denoted by $A=[a_{ij}]\in \mathbb{R}^{n\times n}$, where $a_{ii}=0$ and $a_{ij}\geq 0$ ($a_{ij}>0$ if and only if there is an edge from agent $j$ to agent $i$). The Laplacian $L=[l_{ij}]\in \mathbb{R}^{n\times n}$ of digraph $\mathcal{G}$ is defined as $l_{ii}=\sum_{j\neq i}a_{ij}$ and $l_{ij}=-a_{ij} (j\neq i)$. Define the in-degree and out-degree of node $i$ as $d_i^{in}=\sum_ja_{ij}$ and $d^{out}_i=\sum_ja_{ji}$, respectively. Node $i$ is balanced if and only if its in-degree equals its out-degree, the digraph $\mathcal {G}$ is weight-balanced if and only if all of its nodes are balanced. Note that $L\textbf{1}=0$, a digraph is weight-balanced if and only if $\textbf{1}^\top L=0$. For a digraph, its underlying graph is the graph obtained by replacing all the directed edges with undirected edges. If between every pair of distinct vertices, there is a directed path, this digraph is said to be strongly connected. A digraph is weakly connected if its underlying graph is connected. As stated in \cite{mesbahi2010graph}, let $\mathcal{G}$ be a digraph with Laplacian matrix $L$, then $\hat L=L+L^\top$ is a valid Laplacian matrix for its mirror graph $\hat{\mathcal{G}}$ if and only if $\mathcal{G}$ is balanced. A digraph is said to be undirected if $a_{ij}=a_{ji}$ ($i,j=1,{\dots},n$). Obviously, any undirected graph is balanced.

A dynamic system
\begin{align}\label{sys:def}
  \dot{x}=f(x,u),\quad  y=h(x),\quad x\in \mathbb{R}^n,\, u,y\in\mathbb{R}^p
\end{align}
is said to be passive (\cite{schaft1999l2}), if there exists a continuously
differentiable function $V(x)\geq 0$ such that
\begin{align*}
  \dot{V}\leq -W(x)+u^\top y
\end{align*}
for some positive semidefinite function $W(x)$. $V(x)$ is often called its
storage function. This system is said to be strictly passive if $W(x)$ is
positive definite.  Passivity, due to its explicit physical meaning and
simplicity to manipulate, has been extensively by many authors both for a
single plant and multi-agent systems (e.\,g., \cite{bai2011cooperative,schaft1999l2}).  For simplicity,  we only  consider in this paper the single-input single-output case, i.\,e., $p=1$.

It has been proved that the passivity of a dynamic system is much related to
its high-frequency gain (\cite{khalil2002nonlinear}), which represents the
motion direction of the system in any control strategy.  In most of existing
literatures, this control direction is assumed to be known a prior, or the
high-frequency-gain sign is positive, e.\,g.,\cite{bai2011cooperative,schaft1999l2,su2013cooperative}. As having been mentioned, the control direction of a plant might be unknown or change
under sudden structural damages. Thus, we consider a general class of
nonlinear passive systems with unknown control directions, that is, for
system \eqref{sys:def}, there exists a continuously differentiable function
$V(x)\geq 0$ such that
\begin{align*}
  \dot{V}(x)\leq -W(x)+buy
\end{align*}
where $b$ is an unknown non-zero constant and $W(x)$ some positive semidefinite function.

It can be found that if the constant $b$ or its sign is known a prior, then let $\bar u=\mbox{sign}(b)u$, this system is passive with input $\bar u$ and output $by$. However, if $b$ or its sign is unknown, the conventional passivity-based controller is no longer applicable.

\section{Problem Formulation}

Consider a multi-agent system consisting of $N$ nonlinear agents described by
\begin{align}\label{sys:follower}
  \dot{x}_i=f_i(x_i,u_i),\quad y_i=h_i(x_i),\quad i=1,\dots,N
\end{align}
where $x_i\in \mathbb{R}^{n_i}$,\,$u_i\in\mathbb{R}$,\,$y_i\in \mathbb{R}$
are its state, input, and output of agent $i$. $f_i(\cdot)$ and $h_i(\cdot)$
are locally Lipschitz. Assume that, after a possible inner-loop control,
these agents are all passive but with unknown control directions, i.\,e., for
each $i$, there exist a continuously differentiable function $V_i(x_i)\geq 0$
and a positive semidefinite function $W_i(x_i)$ such that
\begin{align}
  \dot{V}_i(x_i)\leq -W_i(x_i)+b_iu_iy_i
\end{align}
where $b_i$ is an unknown nonzero constant.

Associated with this multi-agent system, a digraph $\mathcal{G}$ can be
defined with the nodes $\mathcal{N}=\{1,...,N\}$ to describe the
communication topology. If the control $u_i$ can get access to the output
information of agent $j$, there is an weighted edge $(j,i)$ in the graph
$\mathcal{G}$, i.\,e., $a_{ij}>0$.

Our control objective is to design $u_i$ for each agent in graph
$\mathcal{G}$, such that output consensus of this multi-agent system composed
of \eqref{sys:follower} can be achieved, i.\,e., $y_i-y_j\to 0$ as $t \to
\infty$ for any $i,\,j=1,\dots,N$ while the overall system maintains bounded.

\begin{remark}\label{rem:formulation:heter}
    In most of existing works considering the coordination of linear
    or nonlinear multi-agent systems \cite{hong2006tracking,ren2008distributed, su2013cooperative, wang2014consensus}, the high-frequency-gain sign of each agent is assumed to be known in prior. But in our formulation, the prior knowledge of each agent's high-frequency-gain sign is no longer necessary. Compared with the results in \cite{chen2014adaptive} where all the high-frequency gains should have an identical sign, agents considered here may have different and unknown control directions.
\end{remark}

\begin{remark}\label{rem:formulation:nonlinear}
    The class of systems considered here can not only strictly cover the systems considered in \cite{peng2014cooperative}, but also include a rather general class of linear or nonlinear passive systems(\cite{schaft1999l2}) with unknown control directions. While passivity has been employed as a powerful tool for group coordination (\cite{bai2011cooperative,chopra2006passivity}), this formulation will enlarge its applications in multi-agent systems.
\end{remark}

\section{Main Results}

In this section, we will first present a Nussbaum-type protocol to achieve leaderless output consensus for those heterogeneous agents, and then provide an extension to the leader-following cases.
\begin{theorem}\label{thm:thm1}
  Consider the multi-agent system consisting of $N$ agents given by \eqref{sys:follower}, there exists a distributed adaptive controller of the form
  \begin{align}\label{eq:ctr11}
    u_i=-\mathcal{N}(k_i)\xi_i, \quad \dot{k}_i=y_i\xi_i
  \end{align}
  where $\xi_i=\sum_{j=1}^N a_{ij}(y_i-y_j)$ and $\mathcal{N}(k_i)=k_i^2\sin(k_i)$,
  such that the output consensus of this multi-agent system is achieved when the communication graph is undirected and strongly connected.
\end{theorem}

{\noindent \bf Proof}:
The proof will be spit into two steps.\medskip

\emph{Step 1:} We first prove the boundedness of $x_i$, $u_i$, and $k_i$.  From the smoothness of related functions, the solution of this closed-loop system will be well-defined on its maximal interval $[0,\,t_f)$. We claim that $t_f=+\infty$. In the following, we will prove it by seeking a contradiction. At first, assume $t_f$ is finite. Taking $V_i(x_i)$ as a sub-Lyapunov function gives
\begin{align}\label{eq:thm1:eq1}
  \dot{V}_i(x_i(t))\leq -W_i(x_i(t))-b_i\mathcal{N}(k_i)y_i\xi_i\leq- b_i\mathcal{N}(k_i)\dot{k}_i.
\end{align}
By integrating both sides from $0$ to $t$, it follows
\begin{align}\label{eq:thm1:eq2}
  {V}_i(x_i(t))-V_i(x_i(0)) \leq -\int_{k_i(0)}^{k_i(t)}b_i\mathcal{N}(s)\, {\rm d}s
\end{align}
and hence
\begin{align}\label{eq:thm1:eq3}
  {V}_i(x_i(t)) \leq  b_i[k_i(t)^2\cos(k_i(t))-2k_i(t)\sin(k_i(t))-2\cos(k_i(t))]+C
\end{align}
where $C=V_i(x_i(0))-b_i[k_i(0)^2\cos(k_i(0))-2k_i(0)\sin(k_i(0))-2\cos(k_i(0))]$ is a finite consonant.

We will prove the boundedness of all trajectories to get a contradiction. To prove this, we now seek another contradiction. Without loss of generality, suppose $k_i(t)$ is upper unbounded. From the continuousness of $k_i(t)$, we can choose a monotonic increasing sequence ${t^i_n}$ such that
\begin{align*}
  k_i(t^i_n)=\begin{cases}
    (2n+1)\pi,\quad \text{if } b_i>0,\\
    2n\pi,\quad \text{if } b_i<0.
  \end{cases}
\end{align*}
Clearly, for any $i$, $\lim_{n\to +\infty}t_n^i=t_f$.

By direct calculations, one has
\begin{align}
V_i(x_i(t_n^i)) \leq \begin{cases}
  b_i[-(2n+1)^2\pi^2+2]+C, \quad \text{if } b_i>0,\\
  b_i[4n^2\pi^2-2]+C,\quad \text{if } b_i<0.
\end{cases}
\end{align}
From this, we can deduce that $V_i(x_i(t_n^i))<0$ for a large enough $n$,
which contradicts the positive semidefiniteness of $V_i(x_i)$.  Therefore,
$k_i(t)$ is bounded during $[0, \,t_f)$ for each $i$. From \eqref{eq:ctr11}
and \eqref{eq:thm1:eq3}, $x_i(t)$, $u_i(t)$ and $\dot{x}_i(t),\,
\dot{k}_i(t)$ are also bounded during $[0, \,t_f)$ for each $i$.  This
implies with a contradiction  argument that no finite-time escape phenomenon
happens and $t_f=\infty$.\medskip

\emph{Step 2:} In this step, we will show that $y_i-y_j\to 0 \,(t \to \infty)$ for any $i,\,j=1,\dots,N$.
Let $K(t)=\sum_{i=1}^N\dot{k}_i(t)$, then
\begin{align*}
  \dot{K}(t)=\sum_{i=1}^N\ddot{k}_i(t)=\sum_{i=1}^N[\dot{y}_i\xi_i+y_i\dot{\xi_i}]
\end{align*}
which is bounded from the boundedness of $y_i(t)$, $\xi_i(t)$ and $\dot{x}_i(t)$. As a result, $K(t)$ is uniformly continuous with respect to time $t$ during $[0,\,\infty)$.

Also, note that
\begin{align}
  \int_0^\infty\,K(t)\,{\rm d}t=\sum_{i=1}^N[{k}_i(\infty)-k_i(0)]\leq K^*
\end{align}
where $K^*$ is a finite constant determined by the bound of $k_i(t)$ ($i=1,\dots,N$). That is, $K(t)$ is integrable on $[0,\,\infty)$. By Barbalat's Lemma, we can derive
\begin{align*}
  \lim_{t\to\infty}{K}(t)=\lim_{t\to\infty}\sum_{i=1}^N\dot{k}_i(t)=\lim_{t\to\infty}y^\top Ly=0
\end{align*}
where $y=\mbox{col}(y_1,\dots,y_N)$ and $L$ is the Laplacian of the
undirected graph. From its connectivity assumption and by Proposition 3.8 in
\cite{mesbahi2010graph}, $y\to \mbox{span}\{ \mathbf{1}_N\}$, which completes
the proof.
\hfill\rule{4pt}{8pt}

\bigskip

In the control law \eqref{eq:ctr11}, multiple Nussbaum gains are employed to tackle the problem of unknown heterogeneous high-frequency-gain signs. Although we have chosen $k^2\sin(k)$ as the Nussbaum-type function, it can be verified that, other Nussbaum-type functions such as $k^2\cos(k)$ and $e^{k^2}\cos(k)$ can be employed in Step 1 as well to achieve the leaderless output consensus among this multi-agent system.

\begin{remark}\label{rem:thm1:heter}
 A similar problem was considered in \cite{chen2014adaptive} for pure integrators with unknown control directions. However, those results heavily relied on the assumption that all agents have the same control direction. Here by constructing sub-Lyapunov functions mentioned above, the control design and analysis for this group of interacting systems are significantly simplified.
\end{remark}

\begin{remark}\label{rem:thm1:nonlinear}
This result provides a sufficient condition to the consensus problem among a general class of agents with heterogeneous control directions, which strictly includes the single integrator in \cite{peng2014cooperative} as a special case. Also, since we only assumed those agents having a passivity-like property, this approach may considerably enlarge the applications of passivity as a design tool (\cite{bai2011cooperative}) in multi-agent systems. Moreover, when all agents are homogenous and share an incremental observability property (\cite{hamadeh2012global}), they will eventually achieve state consensus under this control law, which can be taken as an extension to the results in \cite{hamadeh2012global} for networked nonlinear agents with unknown control directions.
\end{remark}

Unfolding the control law in \eqref{eq:ctr11}, we may extend this conclusion to a general class of digraphs as follows.

\begin{corollary}\label{coro:thm1:digraph}
   Consider the multi-agent system consisting of $N$ agents given by \eqref{sys:follower}, there exists a distributed adaptive controller of the form
  \begin{align}\label{eq:ctr12}
    u_i=-\mathcal{N}(k_i)\xi_i, \quad \dot{k}_i=y_i\xi_i
  \end{align}
  where $\xi_i=\sum_{j=1}^N a_{ij}(y_i-y_j)$ and $\mathcal{N}(k_i)=k_i^2\sin(k_i)$,
  such that the output consensus of this multi-agent system is achieved when the communication digraph is balanced and weakly connected.
\end{corollary}

{\noindent \bf Proof}:
By similar arguments as in Theorem \ref{thm:thm1}, we can derive that $k_i$,
$x_i$, $\xi_i$, $\dot{x}_i$, and $\dot{\xi}_i$ are bounded and
$\dot{K}(t)=\sum_{i=1}^N\ddot{k}_i(t)$ is uniformly continuous with respect
to time $t$ and integrable on $[0,\,\infty)$.  Hence,
\begin{align*}
  \lim_{t\to\infty}{K}(t)=\lim_{t\to\infty}\sum_{i=1}^N\dot{k}_i(t)=\frac{1}{2}\lim_{t\to\infty}y^\top (L+L^\top)y=0.
\end{align*}
Since a weakly connected balanced digraph is automatically strongly
connected, and hence it contains a rooted out-branching. By Theorem 3.12 in
\cite{mesbahi2010graph}, the corresponding protocol can achieve output
consensus of these agents.
\hfill\rule{4pt}{8pt}

\bigskip

Without additional conditions, the controller in Theorem \ref{thm:thm1} can
only ensure the output consensus, i.\,e., $y_i\to y_j$ ($t\to \infty$) for
$i, j$.  In some problems, we may also need to drive all $y_i$ to some
desired value, such as the equilibrium. It is remarkable that this may happen
automatically in some special cases after applying Barbalat's Lemma and then
deliberately checking its asymptotical property of this multi-agent system as
a whole one. To guarantee this property, we can employ a leader-following
formulation and propose a distributed control law, where the reference is
described by a leader. Unlike in the centralized/decentralized cases when
each agent knows this value, only a few agents are assumed to know it in our
cases to save the communication resources.  With some modifications on the
protocol \eqref{eq:ctr11}, a Nussbaum-type adaptive controller will be
designed for each agent in the network to realize this goal.

To keep consistence, we assume as usual the reference point is generated by a leader (denoted as 0) \begin{align}\label{sys:leader}
\dot{x}_0=0,\quad y_0=x_0.
\end{align}
With the help of graph notations, the information flow between the leader and those other agents can be defined as well.

The following theorem shows how this problem can be solved for those nonlinear heterogeneous multi-agent systems without knowing the control directions.

\begin{theorem}\label{thm:thm2}
  Consider the multi-agent system consisting of $N$ followers given by \eqref{sys:follower} and a leader \eqref{sys:leader} with $x_0(0)=0$, there exists a distributed controller of the form
  \begin{align}\label{eq:ctr21}
    u_i=-\mathcal{N}(k_i)\xi_i, \quad \dot{k}_i=y_i\xi_i
  \end{align}
  where $\xi_i=\sum_{j=1}^N a_{ij}(y_i-y_j)+a_{i0}(y_i-y_0)$ and $\mathcal{N}(k_i)=k_i^2\sin(k_i)$,
  such that the consensus problem of this multi-agent system can be solved when the induced subgraph of those followers are undirected, strongly connected and the leader is globally reachable from any other agent.
\end{theorem}

{\noindent \bf Proof}:
The proof is similar with that of Theorem \ref{thm:thm1}. Following those
procedures, we can first prove the boundedness of $x_i(t)$, $\xi_i(t)$,
$\dot{x}_i$, and $\dot{k}_i(t)$, and hence the uniform continuous of
${K}(t)=\sum_{i=1}^N\dot{k}_i(t)$ with respect to the time $t$. By Barbalat's
Lemma, one can obtain
\begin{align*}
  \lim_{t\to\infty}{K}(t)=\frac{1}{2}\lim_{t\to\infty}e^\top(H+H^\top)e=0
\end{align*}
where $e=\mbox{col}(y_1-y_0,\dots, y_N-y_0)=y^\top$, and $H$ is the submatrix
by removing its first row and column of the Laplacian corresponding to the
$N+1$ networked agent systems. By Lemma 3 in \cite{hong2006tracking}, $e$
will eventually vanish. Thus the proof is completed.
\hfill\rule{4pt}{8pt}

\bigskip

This result still holds for a class of digraphs as follows.
\begin{corollary}\label{coro:thm2:digraph}
  Consider the multi-agent system consisting of $N$ followers given by \eqref{sys:follower} and one leader \eqref{sys:leader} with $x_0(0)=0$, there exists a distributed controller of the form
  \begin{align}\label{eq:ctr22}
    u_i=-\mathcal{N}(k_i)\xi_i, \quad \dot{k}_i=y_i\xi_i
  \end{align}
  where $\xi_i=\sum_{j=1}^N a_{ij}(y_i-y_j)+a_{i0}(y_i-y_0)$, and $\mathcal{N}(k_i)=k_i^2\sin(k_i)$,
  such that the consensus problem of this multi-agent system can be solved when the induced subgraph of those followers are balanced, weakly connected and the leader is globally reachable from any other agent.
\end{corollary}
The \ p\,r\,o\,o\,f \ is similar with that of Theorem \ref{thm:thm2} and thus
omitted.

\begin{remark}\label{rem:thm2:nonlinear}
While only leaderless consensus was investigated in \cite{chen2014adaptive}, we provided sufficient conditions to achieve leader-following consensus.  The results in \cite{peng2014cooperative} are its special cases when all considered agents are single integrators. Moreover, if all agents are constant incremental passive (\cite{jayawardhana2005tracking}), it can be extended to the non-equilibrium cases. When no agent has an unknown control direction, these leader-following consensus results are consistent with that in \cite{ren2008distributed}.
\end{remark}

\section{Simulations}
In this section, we proposed several examples to verify the distributed
designs in Section~4.

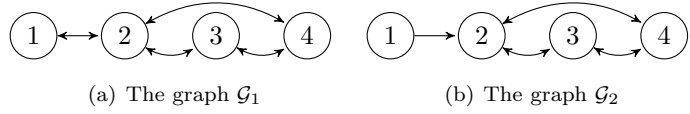
\begin{figure}
  \centering
  \subfigure[The graph $\mathcal{G}_1$]
    {\centering
\begin{tikzpicture}[shorten >=1pt, node distance=1.2 cm, >=stealth',
every state/.style ={circle, minimum width=0.2cm, minimum height=0.2cm}, auto]
\node[align=center,state](node1) {1};
\node[align=center,state](node2)[right of=node1]{2};
\node[align=center,state](node3)[right of=node2]{3};
\node[align=center,state](node4)[right of=node3]{4};
\path[<->]   (node1) edge (node2)
            (node2) edge [bend left] (node4)
            (node2) edge [bend right] (node3)
            (node3) edge [bend right]  (node4)
            ;
\end{tikzpicture}
}\quad
  \subfigure[The graph $\mathcal{G}_2$]
   {\centering
\begin{tikzpicture}[shorten >=1pt, node distance=1.2 cm, >=stealth',
every state/.style ={circle, minimum width=0.2cm, minimum height=0.2cm}, auto]
\node[align=center,state](node1) {1};
\node[align=center,state](node2)[right of=node1]{2};
\node[align=center,state](node3)[right of=node2]{3};
\node[align=center,state](node4)[right of=node3]{4};
\path[<->]   (node1) edge[->] (node2)
            (node2) edge [bend left] (node4)
            (node2) edge [bend right] (node3)
            (node3) edge [bend right]  (node4)
            ;
\end{tikzpicture}
}
\caption{The communication graphs.}\label{fig:graph}\vspace{-5mm}
\end{figure}

\begin{figure}[h!]
  \centering
  \includegraphics[width=0.7\textwidth]{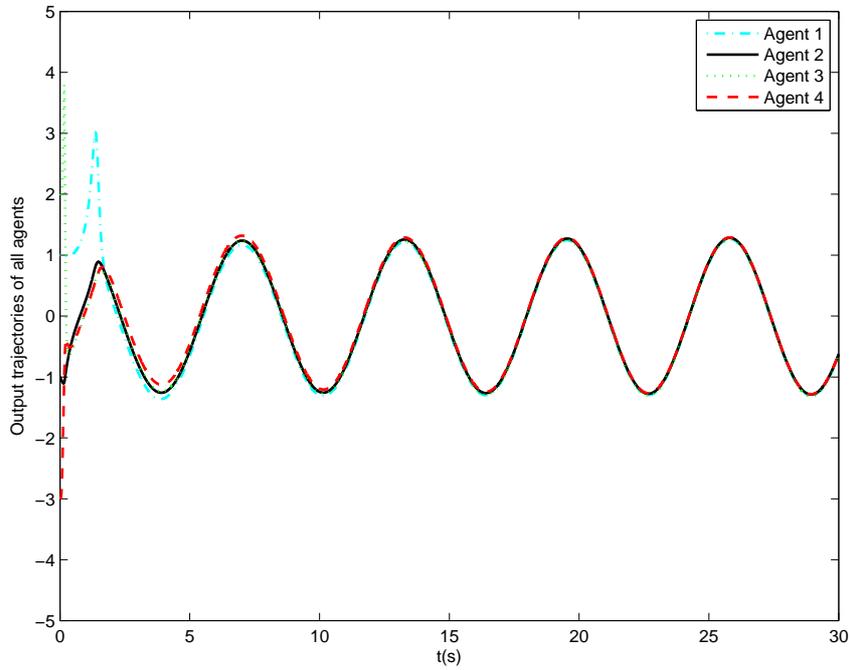}\\ 
  \caption{Output trajectories of agents under control law \eqref{eq:ctr11}.}\label{fig:simu1}
\end{figure}

First, consider a group of controlled oscillators with heterogeneous control directions as follows.
\begin{align*}
\begin{cases}
  \dot{x}_{i1}=x_{i2}\\
  \dot{x}_{i2}=-x_{i1}+b_i u_i\\
y_i=x_{i2}.
\end{cases}
\end{align*}
Their interconnection topology is depicted in Figure~\ref{fig:graph}(a) as an
undirected graph, which satisfies the assumptions in Theorem \ref{thm:thm1}.
Since $b_1, \,b_2, \,b_3,\, b_4$ are unknown and may not have the identical
sign, the control laws in \cite{chen2014adaptive,wang2014consensus} are not
applicable. While it can be verified all these systems share the
passivity-like property, we applied the controller \eqref{eq:ctr11} and its
performance is depicted in Figure~\ref{fig:simu1}.

\begin{figure}[h!]
  \centering
  \includegraphics[width=0.70\textwidth]{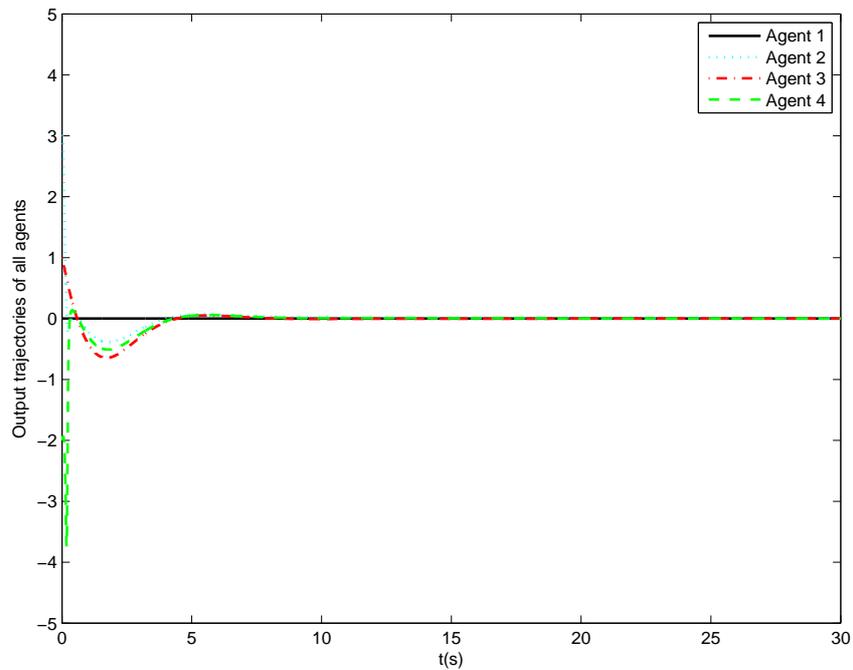}\\ 
  \caption{Output trajectories of agents under control law \eqref{eq:ctr21}.}
  \label{fig:simu2}
\end{figure}

To make it more interesting, we now consider four heterogeneous agents, including a single integrator $\dot{x}_1=0,\, y_1=x_1$ as the leader, two oscillators and a controlled Lorenz system as the followers.
\begin{align*}
\begin{cases}
  \dot{x}_{21}=x_{22}\\
  \dot{x}_{22}=-x_{22}+b_2u_2\\
  y_2=x_{22},
\end{cases}\,
\begin{cases}
\dot{x}_{31}=x_{32}\\
\dot{x}_{32}=-x_{31}-x_{32}+b_3u_3\\
y_3=x_{32},
\end{cases}\,
\begin{cases}
  \dot{x}_{41}=x_{42}-x_{41}\\
  \dot{x}_{42}=x_{41}-x_{42}-x_{41}x_{43}+b_4u_4\\
  \dot{x}_{43}=x_{41}x_{42}-x_{43}\\
  y_4=x_{42}.
  \end{cases}
\end{align*}
Their interconnection topology is depicted in Figure~\ref{fig:graph}(b) by
removing the edge pointed to agent 1, which satisfies the assumptions in
Theorem \ref{thm:thm2}. The simulations result under the controller
\eqref{eq:ctr21} is presented in Figure~\ref{fig:simu2}.
\enlargethispage{5mm}

\section{Conclusions}
An output consensus problem was solved for a general class of nonlinear heterogeneous systems without a prior knowledge of the agents' control directions. Both the leaderless and leader-following consensus were achieved with the help of two distributed Nussbaum-type control laws. Further work will include the extensions for more general systems and graphs with possible disturbances.

\end{document}